\input amstex
\documentstyle{amsppt}
\pagewidth{6.5truein}
\pageheight{9.0truein}
\NoBlackBoxes

\def\cee{{\Bbb C}}

\def\real{{\Bbb R}}
\def\zed{{\Bbb Z}}
\def\cp{{\Bbb C}{\Bbb P}}

\def\Im{\operatorname{Im}}

\def\ep{\varepsilon}
\def\tJ{\widetilde J}
\def\tT{\widetilde T}
\def\tX{\widetilde X}
\def\tomega{\widetilde\omega}
\def\hS{\widehat S}
\def\hX{\widehat X}
\def\homega{\widehat\omega}
\topmatter
\title On symplectically aspherical manifolds with nontrivial $\pi_2$\endtitle
\author Robert E. Gompf\endauthor
\date August 20, 1998\enddate
\address Department of Mathematics, University of Texas at Austin,
Austin TX 78712 \endaddress
\email gompf\@math.utexas.edu \endemail
\thanks Partially supported by NSF grant DMS-9625654. 
\endthanks
\abstract We construct closed symplectic manifolds for which spherical 
classes generate arbitrarily large subspaces in 2-homology, such that 
the first Chern class and cohomology class of the symplectic form both 
vanish on all spherical classes. We construct both K\"ahler and non-K\"ahler 
examples, and show independence of the conditions that these two cohomology 
classes vanish on spherical homology. In particular, we show that the 
symplectic form can pair trivially with all spherical classes even when the 
Chern class pairs nontrivially.
\endabstract 
\endtopmatter

\document

Floer's original approach to the Arnol'd Conjecture \cite{F} worked for 
closed symplectic manifolds $(X,\omega)$ satisfying the additional 
hypotheses that $[\omega]$ and $c_1(\omega)$ in $H^2(X;\real)$ should 
vanish on all spherical homology classes. 
While more recent work has avoided these hypotheses in this particular 
context, elsewhere the 
question has persisted of which symplectic manifolds satisfy the above 
hypotheses. 
Clearly the hypotheses are true for aspherical manifolds, but there are 
apparently no other published examples for which $[\omega]$ vanishes 
on all spherical classes, prompting various people to ask whether such 
``symplectically aspherical'' manifolds are always aspherical. 
In this note, we answer the question in the negative. 
(The author has been informed that J.~Koll\'ar has also 
reached this conclusion --- see Remark~8.) 
To be more precise, let $\Pi (X) = \text{span}_{\real}(\Im (\pi_2 (X) 
\to H_2(X;\real)))$ denote the subspace of $H_2(X;\real)$ 
spanned by spherical classes. 
We construct families of closed, symplectic 4-manifolds $(X,\omega)$ for which 
$\Pi(X)$ has arbitrarily large dimension, but both $\omega$ and $c_1(\omega)$ 
vanish on $\Pi(X)$. 
We produce such a family for which each $\omega$ is a K\"ahler form, and 
a different family for which no $X$ is homotopy equivalent to a K\"ahler 
manifold. 
We also show that the conditions that $\omega$ and $c_1(\omega)$ vanish on 
$\Pi(X)$ are independent --- in particular, we can have $\omega|\Pi(X)=0$ 
but $c_1(\omega)|\Pi(X) \ne0$. 
Note that since $\Pi(X) = H_2(X;\real)$ when $X$ is simply connected, 
any closed symplectic 
manifold with $\omega|\Pi(X)=0$ must have infinite fundamental group. 
The question of which topologies (e.g. fundamental groups) are realized by 
symplectically aspherical manifolds seems highly nontrivial. 

Our construction is based on symplectic branched coverings. 
As observed by Gromov (\cite{G}, 3.4.4(E)), 
symplectic structures can be lifted under 
branched coverings with symplectic branch loci. 
The following lemma pins down the details. 
We denote the Poincar\'e duality isomorphism by $PD$. 

\proclaim{Lemma 1} 
Let $(X,\omega)$ be a closed symplectic manifold with a codimension-2 
symplectic submanifold $B$. 
Suppose that $\pi:\tX\to X$ is a branched 
covering map with branch locus $B$, and  
the manifold $\pi^{-1}(B)
\subset \tX$  has components $B_1,\ldots,B_n$, oriented via the 
symplectic structure on $B$. 
Let $d_i>0$ denote the multiplicity of $\pi$ along $B_i$. 
Then there is a symplectic form $\tomega$ on $\tX$ with 
$[\tomega] = \pi^* [\omega] \in H^2 (\tX;\real)$ and 
$c_1 (\tomega) = \pi^* c_1(\omega) + \sum (1-d_i) PD[B_i]\in H^2 (\tX;\zed)$. 
\endproclaim 

\demo{Proof} 
The form $\pi^*\omega$ is symplectic on $\tX-\bigcup_{d_i\ge2} B_i$. 
By standard symplectic topology (\cite{MS} Example~5.10 and 
Theorem~3.29), we can find a closed tubular 
neighborhood $N$ of $B$ in $X$ with a Hamiltonian circle action, whose  
fixed set is $B= H^{-1}(0)$ for the Hamiltonian $H:N\to [0,\ep]$. 
Symplectic reduction of each $H^{-1}(t)$ gives a family $\omega_t$, 
$t\in [0,\ep]$, of symplectic forms on $B$, with $\omega_0 = \omega|B$. 
For a suitable family $\alpha_t$ of 1-forms on the total space of 
the circle bundle $p:\partial N
\to B$, the formula $p^*\omega_t + \alpha_t\wedge dt$ defines a symplectic form 
on $[0,\ep]\times \partial N$, and symplectic cutting at $t=0$ recovers 
the original form $\omega$ on $N$ (\cite{L}, cf. also \cite{MS} Example~5.10). 
Now we use $\pi$ to lift this structure to a tubular neighborhood $N_i$ 
of each component $B_i$ of $\pi^{-1}(B)$. 
The form and $S^1$-action on $[0,\ep]\times \partial N$ lift to a 
form $\omega'$ and Hamiltonian circle action on $[0,\ep]\times\partial N_i$ 
with Hamiltonian $d_iH\circ \pi$. 
Since $\omega'$ agrees with $\pi^*\omega$ on $(0,\ep]\times \partial N_i 
\approx N_i-B_i$, symplectic cutting at $0$ yields a symplectic form 
$\tomega$ on $\tX$ that agrees with $\pi^*\omega$ outside of $\bigcup N_i$ 
(and with $\tomega = \pi^*(\omega|B)$ on each $B_i$).  
To check that $[\tomega] = \pi^*[\omega]$, we evaluate both on an 
arbitrary homology class, which we can assume is represented by a surface 
$F\subset \tX$ intersecting each $N_i$ in a union of $S^1$-invariant normal 
disks. 
For each such disk $D$, it is easy to see that 
$\int_D\tomega = d_i \int_{\pi(D)}\omega$, so since $\tomega = \pi^*\omega$ 
on $\tX - \bigcup \text{ int }N_i$, we have $\langle \tomega,[F]\rangle = 
\langle \omega,\pi_* [F]\rangle = \langle \pi^* \omega,[F]\rangle$ 
as required. 

To compute $c_1(\tomega)$, fix an almost-complex structure $J$ on $X$ 
tamed by $\omega$, with $B$ and the fibers of $N$ $J$-holomorphic. 
Then $J$ lifts to a $C^0$ almost-complex structure $\tJ$ on $\tX$ tamed by 
$\tomega$. 
Now $c_1(\omega) = c_1(J)$ is Poincar\'e dual to the zero locus of a generic 
section of $\Lambda TX$, where $\Lambda$ denotes 
the complex top exterior power. 
On $N$, we have $TX = TB\oplus \nu B$ and $\Lambda TX = \Lambda TB\otimes 
\nu B$. 
Let $\sigma :B\to \Lambda TB$ be a generic section, and let $\tau :N\to\nu B$ 
be given by the Tubular Neighborhood Theorem (so $\tau^{-1}(0)=B)$. 
Extend the section $\sigma\otimes\tau$ of $\Lambda TX|N$ generically over 
all of $X$. 
The resulting zero locus exhibits $c_1(J)$ as $PD[C] + PD[B]$, for some 
codimension-2 submanifold $C$ of $X$ intersecting $B$ transversely. 
Since $\tau$ lifts to a continuous map $N_i\to \nu B_i$ vanishing 
with multiplicity~1 on $B_i$, the corresponding section of $\Lambda T\tX$ 
exhibits $c_1(\tomega) = c_1(\tJ)$ as 
$\pi^* PD[C] + \sum PD[B_i]  
= \pi^* (c_1(J) - PD[B]) + \sum PD[B_i] 
= \pi^* c_1 (J) + \sum (1-d_i) PD[B_i]$.\qed
\enddemo 

\example{Example 2} 
Let $X= F_1\times F_2$ be a product of two surfaces of nonzero genus, with a 
product symplectic form $\omega$. 
Fix integers $m_1,m_2 \ge1$ and $d\ge2$. 
Choose $m_1d$ surfaces in $X$ of the form $\{p\} \times F_2$ and $m_2d$ 
surfaces $F_1\times \{q\}$, and let $B^*$ denote their union. 
Let $B\subset X$ be the smooth symplectic surface obtained by smoothing 
the $m_1m_2d^2$ positive double points of $B^*$ via the local model 
$\{z_1z_2=0\}\mapsto \{z_1z_2=\ep\}$. 
Let $\pi :\tX\to X$ be the $d$-fold cyclic branched covering of $X$, branched 
along $B$ with multiplicity $d$ on $\pi^{-1}(B)$. 
(There is only one such covering, up to diffeomorphisms of $X$.)
The following two theorems immediately show that $\dim \Pi (\tX) \ge 
m_1m_2 d^2(d-1)$ but the form $\tomega$ given by Lemma~1 satisfies 
$\tomega |\Pi (\tX) = c_1 (\tomega) |\Pi (\tX) =0$. 
\endexample 

\proclaim{Theorem 3} 
For $(\tX,\tomega)$ as in Lemma~1, if $\omega|\Pi(X)=0$ (e.g., if 
$\pi_2(X) =0$) then $\tomega |\Pi(\tX)=0$. 
If $\pi_2 (X) =0 $ and $\pi^{-1}(B)$ is connected then 
$c_1(\tomega)|\Pi (\tX) =0$. 
\endproclaim 

\demo{Proof} 
If $\alpha \in H_2 (\tX;\zed)$ is represented by a map of a sphere then so 
is $\pi_* \alpha \in H_2 (X;\zed)$. 
Thus $\langle[\tomega],\alpha\rangle = \langle [\omega],\pi_*\alpha\rangle=0$. 
If, in addition, $\pi_2(X)=0$, then $\pi_*\alpha =0$ so the intersection 
number $\pi_* \alpha\cdot B$ vanishes. 
Connectedness of $\pi^{-1}(B)$ now implies that $\alpha\cdot \pi^{-1}(B)=0$ 
and $\langle c_1(\tomega) ,\alpha\rangle = \langle c_1(\omega),\pi_*
\alpha\rangle =0$.\qed 
\enddemo 

\proclaim{Theorem 4} 
Let $\pi:\tX\to X$ be a $d$-fold branched covering of orientable 
4-manifolds, whose branch locus $B\subset X$ is obtained from a 
generically immersed surface $B^*\subset X$ by smoothing all double 
points as in Example~2 (up to orientation). 
If $B^*$ has $k\ge1$ double points then $\Pi (\tX) \ne0$. 
If $\pi$ is injective on $\pi^{-1}(B)$ then $\dim \Pi(\tX) \ge k(d-1)$. 
\endproclaim

\demo{Proof} 
At each double point of $B^*$, a 4-ball neighborhood $K$ intersects $B$ 
in the complex curve $z_1z_2=\ep$, which is the Milnor fiber $M(2,2)$. 
Each component of $\pi^{-1}(K)$ intersects one $B_i$ and is the $d_i$-fold 
cover of $K$ branched along $K\cap B$. 
This is diffeomorphic to the Milnor fiber $M(2,2,d_i)$, which is 
the plumbing along a linear graph of $d_i-1$ spheres of square $-2$. 
The intersection matrix of these $d_i-1$ spheres has nonzero determinant, 
so each component of each $\pi^{-1}(K)$ contributes $d_i-1$ spheres to a 
linearly independent subset of $\Pi (\tX)$. 
If $\pi |\pi^{-1} (B)$ is injective then each $d_i=d$, so we obtain 
$\dim \Pi (\tX) \ge k(d-1)$ (and in general, $\dim \Pi(\tX) \ge kd-\ell$, 
where $\ell$ is the number of components of $\bigcup \pi^{-1}(K)$).\qed
\enddemo

\example{Example 5 -- K\"ahler manifolds} 
We find a family as in Example~2 with $\Pi(\tX)$ arbitrarily large, for 
which each $\tomega$ is a K\"ahler form. 
Suppose that each $F_i$ is given as a complex curve of degree $m_i$ in  
$\cp^{N_i}$.  
Let $H_i \subset \cp^{N_i}$ be a generic degree-$d$ hypersurface, and 
choose the $m_1d$ points $p$ in Example~2 to be the set $F_1\cap H_1$. 
Similarly, choose the points $q$ to be $F_2\cap H_2$. 
Thus, each of these finite subsets of $F_i$ is exhibited as the zero locus of a 
holomorphic section of $L_i^{\otimes d}$, where $L_i\to F_i$ is the 
restriction of the hyperplane bundle of $\cp^{N_i}$. 
The corresponding zero locus of $L_1^{\otimes d} \otimes L_2^{\otimes d} 
\to X= F_1\times F_2$ is the singular surface $B^*$. 
By  varying our choices of $H_i$, it is easy to find three holomorphic 
sections of $L_1^{\otimes d}\otimes L_2^{\otimes d}$ whose common zero 
locus is empty --- thus, the complete linear system of holomorphic sections 
of $L_1^{\otimes d}\otimes L_2^{\otimes d}$ has empty base locus. 
By Bertini's Theorem, a generic holomorphic section $\rho$ of 
$L_1^{\otimes d} \otimes L_2^{\otimes d}$ has smooth zero locus, so the 
required perturbation of $B^*$ to a smooth surface $B= \rho^{-1}(0)$ can be 
done holomorphically. 
To find a K\"ahler form $\tomega$ on the resulting $\tX$ satisfying the 
conclusion of Lemma~1, note that the $d$-fold tensor power map 
$p:L_1\otimes L_2 \to (L_1\otimes L_2)^d = L_1^{\otimes d}\otimes 
L_2^{\otimes d}$ is a $d$-fold covering map branched along the zero section. 
Since the section $\rho :X\to L_1^{\otimes d}\otimes L_2^{\otimes d}$ is 
a diffeomorphism onto its image,  we can identify $\tX$ with 
$p^{-1}(\Im \rho)\subset L_1\otimes L_2$ and $\pi :\tX \to X$ with the 
restriction of projection $L_1\otimes L_2 \to X$. 
If we identify $X$ with the zero section of $L_1\otimes L_2$, then any 
K\"ahler form $\Omega$ on the total space of $L_1\otimes L_2$ restricts 
to K\"ahler forms $\tomega$ on $\tX$ and $\omega$ on $X$. 
For any $\alpha\in H_2(\tX;\real)$, the classes $\alpha$ and $\pi_*\alpha$ 
have the same image in $H_2(L_1\otimes L_2;\real)$, so 
$\langle \tomega,\alpha\rangle = \langle\Omega,\pi_*\alpha\rangle 
= \langle \omega,\pi_* \alpha\rangle$, implying that $[\tomega] = 
\pi^* [\omega]$ as required. 
The Chern class formula follows as before (and is well-known since $\pi$ 
is holomorphic), so we again have 
$\tomega |\Pi (\tX) = c_1(\tomega) | \Pi (\tX) =0$. 
\endexample

\example{Example 6 -- Non-K\"ahler manifolds} 
Let $(X,\omega)$ denote the Kodaira-Thurston manifold, obtained as the 
quotient of $(\real^4,dx_1\wedge dx_2 + dx_3\wedge dx_4)$ by the discrete 
group of symplectomorphisms generated by unit translations in the $x_1,x_2$ 
and  $x_3$ directions together with $(x_1,x_2,x_3,x_4) \mapsto (x_1+x_2, 
x_2,x_3,x_4+1)$. 
The projection $(x_1,x_2,x_3,x_4) \mapsto (x_3,x_4)$ induces a bundle 
structure $p:X\to T^2$ with torus fibers and a section $x_1=x_2=0$ with 
trivial normal bundle. 
As in Example~2, choose $m_1d$ fibers and $m_2d$ parallel copies of the 
section, smooth the intersections of the resulting $B^*$ to obtain 
a smooth symplectic surface $B\subset X$, and let $\pi :\tX \to X$ be the 
resulting $d$-fold branched cover. 
As before, we obtain examples with $\Pi (\tX)$ arbitrarily large, but with 
the form $\tomega$ of Lemma~1 satisfying $\tomega |\Pi (\tX) = c_1(\tomega) 
|\Pi (\tX) =0$. 
This time, however, no $\tX$ has the homotopy type of a K\"ahler manifold,  
as we verify by showing that $b_1(\tX) =3$. 
Recall that $b_1 (X)=3$, since $H_1 (X;\zed)$ has 4 canonical generators 
$\alpha_1,\ldots,\alpha_4$ with $\alpha_i$ descending from the 
$x_i$-axis, and the monodromy introduces 
the relation $\alpha_2 =\alpha_1 +\alpha_2$. 
Since $\pi_* :H_1(\tX;\zed) \to H_1(X;\zed)$ is surjective, we have 
$b_1 (\tX) \ge 3$. 
The Milnor fibers $M(2,2,d) \subset \tX$ arising as in the proof 
of Theorem~4 are simply connected, so we can collapse these to points, 
obtaining a singular space $\tX^*$ with $b_1(\tX^*) = b_1(\tX)$ and a 
map $\pi^* :\tX^* \to X$ that is a $d$-fold covering branched along $B^*$. 
The composite $p \circ \pi^* :\tX^* \to T^2$ is a singular fibration whose 
$m_1d$ singular fibers are tori.  
By transversality, $H_1(\tX^*;\real)$ is generated by cycles disjoint 
from the singular fibers, hence, by cycles in a regular fiber together 
with a pair of classes lifted from the base $T^2$. 
(Meridians to the critical values in $T^2$ lift to nullhomologous 
$\real$-cycles in $\tX^*$.) 
Since any cycle in a regular fiber is homologous to one in a preassigned 
singular fiber, $H_1(\tX^*;\real)$ is generated by 4 classes 
$\tilde\alpha_1,\ldots,\tilde\alpha_4$ lifted from cycles in $B^*$ 
representing $\alpha_1,\ldots,\alpha_4$. 
Since the relation $\alpha_2 = \alpha_1 + \alpha_2$ is represented by a 
surface in $X$ disjoint from $B^*$ and lifting trivially to $\tX^*$, 
we obtain the relation $\widetilde\alpha_2 = \widetilde\alpha_1
+\widetilde \alpha_2$ 
in $H_1(\tX^*;\real)$, so $b_1 (\tX) = b_1 (\tX^*) = 3$ as required. 
\endexample

\proclaim{Theorem 7} 
The conditions $\omega | \Pi (X) =0$ and $c_1(\omega) |\Pi(X)=0$ are 
independent. 
\endproclaim 

\demo{Proof} 
We have shown that both  $\omega$ and $c_1(\omega)$ can vanish on a 
nontrivial $\Pi(X)$. 
A small perturbation of $\omega$ gives examples with $\omega |\Pi (X)\ne0$, 
$c_1(\omega) |\Pi(X)=0$. 
(The $K3$-surface is another such example.) 
Any simply-connected K\"ahler surface other than $K3$ satisfies 
$\omega,c_1(\omega)\ne0$ on $\Pi (X) = H_2(X;\real)$. 
For the remaining case $\omega|\Pi (X)=0$, $c_1(\omega) |\Pi (X)\ne0$, begin 
with the torus $T^4= \cee^2/\zed^4$ with the standard K\"ahler form $\omega$. 
Let $B^* \subset T^4$ be the union of tori descending from the 4 lines 
$z_i = \pm\frac14$ in $\cee^2$, and let $B$ be the symplectic surface 
obtained by smoothing the 4 intersections $(\pm\frac14,\pm\frac14)$ of $B^*$ 
using the local model $z_1z_2 = -\ep$, $0<\ep \ll1$ (up to translation 
and perturbation away from the double points). 
The torus $T\subset T^4$ determined by the equation $z_2 = \bar z_1$ 
intersects $B^*$ in the 2 points $\pm (\frac14,\frac14)$ and is disjoint 
from $B$. 
Since  $T$ is Lagrangian, it is easy to perturb $\omega$ so that $T$  
becomes symplectic.  
Let $(\tX,\tomega)$ be a 2-fold cover of $T^4$ branched along $B$ 
(and chosen to be trivial over $T$), and let $\tT$ be one of the 2 lifts 
of $T$ to $\tX$. 
The vanishing cycle of $B$ created by smoothing $B^*$ at $(\frac14,\frac14)$ 
bounds a disk $D$ in $T^4$ intersecting $T$ transversely in one point. 
(The coordinate change $z_1 = u_1 + iu_2$, $z_2 = -u_1 + iu_2$ (following 
the appropriate translation) exhibits $B$ as $u_1^2 + u_2^2 = \ep$, 
$D$ as the $\sqrt{\ep}$-disk in $\real^2$, and $T$ as  $i\real^2$.) 
The 2 lifts of $D$ to $\tX$ fit together to form a sphere $S\subset \tX$
intersecting $\tT$ transversely in a single point. 
Let $(\hX,\homega)$ be any branched cover of $\tX$ with branch locus given by 
2 parallel copies of $\tT$. 
(Such a branched cover exists since $\tT$ is symplectic with 
trivial normal bundle --- e.g., 
send the two meridians to opposite generators of $\zed/d$.) 
The sphere $S\subset \tX$ lifts to a sphere $\hS \subset \hX$ intersecting 
two components $B_i$ of the lifted branch locus of $\hX\to \tX$ each 
in one point. 
The two components $B_i$ have the same multiplicity $d\ge2$, 
and the signs $B_i\cdot \hS$ agree if each $B_i$ is symplectically oriented. 
Now Theorem~3 implies that $\homega |\Pi (\hX) =0$ and 
$c_1(\tomega) |\Pi (\tX) =0$, so by Lemma~1 applied to $\hX\to \tX$, we have  
$\langle c_1(\homega),\hS\,\rangle = \sum (1-d_i)B_i\cdot\hS = 2(1-d)\ne0$.\qed
\enddemo 

\example{Remark 8} 
J. Koll\'ar produced a different construction of symplectically aspherical 
K\"ahler manifolds with $\Pi(X) \ne 0$, 
which preceded the author's work by several months. 
He now observes that both constructions are based on a common idea: 
If $(X,\omega)$ is symplectic, $f:X\to Y$ is a continuous map to a 
topological space with $\pi_2(Y)=0$, and $[\omega] = f^*\Omega$ for some 
$\Omega\in H^2(Y;\real)$, then $\omega|\Pi(X)=0$. 
(The proof is the same as for Theorem~3.) 
In our case, $f$ is a branched covering, whereas Koll\'ar took $f$ to be 
an inclusion of algebraic varieties with $Y$ a $K(\pi,1)$ (e.g., an 
Abelian variety) with $\dim_{\cee}Y\ge 3$. 
\endexample

\enddocument